\documentclass{article}

\usepackage{amssymb}
\usepackage{amsmath}
\usepackage{tensor}
\newcommand{\proof}{{\bf Proof:  }}
\newcommand{\remark}{{\bf Remark:  }}

\newcommand{\example}{{\bf Example:  }}

\newcommand{\dimv}{\underline{\dim}}

\newcommand{\hb}{\newline\hspace*{\fill}$\Box$}

\newcommand{\PrR}{{}^\prime\!R}
\newcommand{\PrL}{{}^\prime\!\Lambda}
\newtheorem{theorem}{Theorem}[section]
\newtheorem{lemma}[theorem]{Lemma}

\newtheorem{proposition}[theorem]{Proposition}
\newtheorem{corollary}[theorem]{Corollary}

\begin{document}
\parindent0pt
\title{\bf Cohomology of quiver moduli, functional equations, and integrality of Donaldson-Thomas type invariants}

\author{Markus Reineke\\ Fachbereich C - Mathematik\\ Bergische Universit\"at Wuppertal\\ D - 42097 Wuppertal, Germany\\
e-mail: reineke@math.uni-wuppertal.de}
\date{}
\maketitle

\begin{abstract} A system of functional equations relating the Euler characteristics of moduli spaces of stable representations of quivers and the Euler characteristics of (Hilbert scheme-type) framed versions of quiver moduli is derived. This is applied to wall-crossing formulas for the Donaldson-Thomas type invariants of M. Kontsevich and Y. Soibelman, in particular confirming their integrality.
\end{abstract}

\section{Introduction}

In \cite{KS}, a framework for the definition of Donaldson-Thomas type invariants for Calabi-Yau categories endowed with a stability structure is developed. One of the key features of this setup is a wall-crossing formula for these invariants, describing their behaviour under a change of stability structure in terms of a factorization formula for automorphisms of certain Poisson algebras defined using the Euler form of the category.\\[1ex]
In \cite{RWC}, such factorization formulas are interpreted using quiver representations, their moduli spaces, and Hall algebras. The main result of \cite{RWC} interprets the factorization formula in terms of generating series of the Euler characteristic of the smooth models of \cite{SM}, which can be viewed as Hilbert schemes in the setup of quiver moduli:\\[1ex]
In the general framework of \cite{KR,LBNCC}, series of moduli spaces of stable representations of quivers are viewed as the commutative `approximations' to a fictitious noncommutative geometry of (the path algebras of) quivers. In this framework, the smooth models can be viewed as Hilbert schemes of points of this noncommutative geometry (for example, in the case of moduli spaces of semisimple representations of quivers, the smooth models parametrize finite codimensional left ideals in the path algebra of the quiver, in the same way as the Hilbert schemes of points of an affine variety parametrize finite codimensional ideals in the coordinate ring of the variety; see \cite[Section 6]{SM}). Since path algebras of quivers are of global dimension $1$, this setup thus describes aspects of a one-dimensional noncommutative geometry.\\[1ex]
The first aim of this paper (after reviewing some facts on quiver moduli in Section \ref{recoll}) is to develop a (one-dimensional, noncommutative) analog of the result \cite{Ch} calculating the generating series of Euler characteristics of Hilbert schemes of points of a threefold $X$ as the $\chi(X)$-th power of the MacMahon series (see \cite[Theorem 4.12]{BF}, \cite[Conjecture 1]{MNOP} for the corresponding statement for Donaldson-Thomas invariants). Namely, we relate the (generating series of) Euler characteristics of moduli spaces of stable quiver representations and  Euler characteristics of their smooth models by a coupled system of functional equations, see Theorem \ref{t42}, Corollary \ref{corsd}. This is achieved using a detailed analysis of a Hilbert-Chow type morphism from a smooth model to a moduli space of semistable representations, whose fibres are non-commutative Hilbert schemes (see Section \ref{fe1}). The explicit cell decompositions for the latter, constructed in \cite{SM}, yield functional equations for the Euler characteristic; see Section \ref{section4}.\\[1ex]
The second aim is to prove the integrality conjecture \cite[Conjecture 1]{KS} for the Donaldson-Thomas type invariants appearing in the wall-crossing formula of \cite{KS}; see Section \ref{app}. These numbers arise by a factorization of the generating series of Euler characteristics as an Euler product (this process can thus be interpreted as fitting a genuinely noncommutative (one-dimensional) object into a commutative (three-dimensional) framework). Using the functional equations mentioned above, we can interprete this process as passing to the compositional inverse of an Euler product, and elementary number-theoretic considerations in Section \ref{number} yield the desired integrality property (it should be noted that a similar process appears in \cite{Sti} in relating modular forms and instanton expansions). We also confirm a conjectural formula of \cite{KS} for diagonal Donaldson-Thomas type invariants using recent results of \cite{Weist}.\\[2ex]
{\bf Acknowledgments:} I would like to thank T. Bridgeland, V. Jovovic, S. Mozgovoy, Y. Soibelman, H. Thomas, V. Toledano-Laredo and T. Weist for interesting discussions concerning this work.

\section{Recollections on quiver moduli}\label{recoll}

In this section, we fix some notation and collect information on moduli spaces of stable representations of quivers and some of their variants, like Hilbert schemes of path algebras and the smooth models of \cite{SM}. See \cite{Rmoduli} for an overview over these moduli spaces and the techniques used to prove some of the results cited below.\\[1ex]
Let $Q$ be a finite quiver, with set of vertices $I$, and arrows written as $\alpha:i\rightarrow j$ for $i,j\in I$. Denote by $r_{i,j}$ the number of arrows from $i\in I$ to $j\in I$ in $Q$. Define $\Lambda={\bf Z}I$, with elements written in the form $d=\sum_{i\in I}d_ii$, and define $\Lambda^+={\bf N}I\subset \Lambda$. We will sometimes use locally finite quiver, for which the set of vertices is possibly infinite, but with only finitely many arrows starting or ending in each single vertex. Dimension vectors for locally finite quivers are assumed to be supported on a finite subquiver.\\[1ex]
Introduce a non-symmetric bilinear form $\langle\_,\_\rangle$ (the Euler form) on $\Lambda$ by
$$\langle d,e\rangle=\sum_{i\in I}d_ie_i-\sum_{\alpha:i\rightarrow j}d_ie_j$$
for $d,e\in\Lambda$; we thus have $\langle i,j\rangle=\delta_{i,j}-r_{i,j}$. For a functional $\Theta\in\Lambda^*={\rm Hom}_{\bf Z}(\Lambda,{\bf Z})$ (called a stability), define the slope of $d\in\Lambda^+\setminus 0$ as $\mu(d)=\Theta(d)/\dim d$, where $\dim d=\sum_{i\in I}d_i$. For $\mu\in{\bf Q}$, define
$$\Lambda^+_\mu=\{d\in\Lambda^+\setminus 0\, ,\, \mu(d)=\mu\}\cup\{0\}$$
(a subsemigroup of $\Lambda^+$), and $\PrL_\mu^+=\Lambda^+_\mu\setminus 0$.\\[1ex]
We consider complex finite dimensional representations $M$ of $Q$, consisting of a tuple of complex vector spaces $M_i$ for $i\in I$ and a tuple of ${\bf C}$-linear maps $M_\alpha:M_i\rightarrow M_j$ indexed by the arrows $\alpha:i\rightarrow j$ of $Q$. The dimension vector $\dimv M\in\Lambda^+$ is defined by $(\dimv M)_i=\dim_{\bf C}M_i$. The abelian ${\bf C}$-linear category of all such representations is denoted by ${\rm mod}_{\bf C}Q$.\\[1ex]
Define the slope of a non-zero representation $M$ of $Q$ as the slope of its dimension vector, thus $\mu(M)=\mu(\dimv M)$. Call $M$ semistable (for the choice of stability $\Theta$) if $\mu(U)\leq \mu(M)$ for all non-zero subrepresentations $U$ of $M$, and call $M$ stable if $\mu(U)<\mu(M)$ for all proper non-zero subrepresentations $U$ of $M$. Finally, call $M$ polystable if it isomorphic to a direct sum of stable representations of the same slope. The full subcategory ${\rm mod}_{\bf C}^\mu Q$ of all semistable representations of slope $\mu\in{\bf Q}$ is an abelian subcategory, that is, it is closed under extensions, kernels and cokernels. Its simple (resp.~semisimple) objects are precisely the stable (resp.~polystable) representations of $Q$ of slope $\mu$.\\[1ex]
Note that in the case $\Theta=0$, all representations are semistable, and the stable (resp.~polystable) ones are just the simples (resp.~semisimples).\\[2ex]
By \cite{King}, for every $d\in\Lambda^+$, there exists a (typically singular) complex variety $M_d^{\rm sst}(Q)$ whose points pa\-ra\-me\-trize the isomorphism classes of polystable representations of $Q$ of dimension vector $d$. In case $\Theta=0$, the variety $M_d^{sst}(Q)$ is affine, parametrizing isomorphism classes of semisimple representations of $Q$ of dimension vector $d$; it will be denoted by $M_d^{ssimp}(Q)$. This variety always contains a special point $0$ corresponding to the semisimple representations $\bigoplus_{i\in I}S_i^{d_i}$, where $S_i$ denotes the one-dimensional representation of $Q$ concentrated at a vertex $i\in I$, and with all arrows represented by zero maps. Note that all $M_d^{ssimp}(Q)$ reduce to the single point $0$ if $Q$ has no oriented cycles. There exists a projective morphism from $M_d^{sst}(Q)$ to $M_d^{ssimp}(Q)$.\\[1ex]
The variety $M_d^{sst}(Q)$ admits the following Luna type stratification (that is, a finite decomposition into locally closed subsets) induced by the decomposition types of polystable representations: let $\xi=((d^1,\ldots,d^s),(m_1,\ldots,m_s))$ be a pair consisting of a tuple of dimension vectors of the same slope as $d$ and a tuple of non-negative integers, such that $d=\sum_{i=1}^sm_id^i$. We call such $\xi$ a polystable type for $d$. Analogously to \cite{LBP} in the case of trivial stability, the set of all polystable representations $M$ such that $M=\bigoplus_{i=1}^sU_i^{m_i}$ for pairwise non-isomorphic stable representations $U_i$ of dimension vector $d^i$ forms a locally closed subset of $M_d^{sst}(Q)$, denoted by $S_\xi$.\\[1ex]
Let $n\in\Lambda^+$ be another dimension vector, and fix complex vector spaces $V_i$ of dimension $n_i$ for $i\in I$. A pair $(M,f)$ consisting of a semistable representation $M$ of $Q$ of dimension vector $d$ and a tuple $f=(f_i:V_i\rightarrow M_i)$ of ${\bf C}$-linear maps is called stable in \cite{SM} if the following condition holds: if $U$ is a proper subrepresentation of $M$ containing the image of $f$ (in the sense that $f_i(V_i)\subset U_i$ for all $i\in I$), then $\mu(U)<\mu(M)$. Two such pairs $(M,f)$, $(M',f')$ are called equivalent if there exists an isomorphism $\varphi:M\rightarrow M'$ interwining the additional maps, that is, such that $f_i'=\varphi_i\circ f_i$ for all $i\in I$.\\[1ex]
By \cite{SM}, there exists a smooth complex variety $M_{d,n}^\Theta(Q)$, called a smooth model for $M_d^{sst}(Q)$, whose points parametrize equivalence classes of stable pairs as above. It admits a projective morphism $\pi_d:M_{d,n}^\Theta(Q)\rightarrow M_d^{sst}(Q)$.\\[1ex]
In the case of trivial stability, the smooth model (a Hilbert scheme for the path algebra of $Q$) ${\rm Hilb}_{d,n}(Q):=M_{d,n}^0(Q)$ parametrizes arbitrary representations $M$ of $Q$ of dimension vector $d$, together with maps $f_i:V_i\rightarrow M_i$ whose images generate the representation $M$. There exists a projective morphism $\pi:{\rm Hilb}_{d,n}(Q)\rightarrow M_d^{ssimp}(Q)$. We denote by ${\rm Hilb}_{d,n}^{nilp}(Q)$ the inverse image under $\pi$ of the special point $0\in M_d^{ssimp}(Q)$; it parametrizes pairs $(M,f)$ as above, with $M$ being a nilpotent representation, in the sense that all maps $M_{\alpha_n}\circ\ldots\circ M_{\alpha_1}$ representing oriented cycles $i_1\stackrel{\alpha_1}{\rightarrow}i_2\stackrel{\alpha_2}{\rightarrow}\ldots\stackrel{\alpha_n}{\rightarrow}i_1$ in $Q$ are nilpotent.\\[1ex]
Following \cite{ALB}, for any polystable type $\xi$ for $d$ as above, introduce new (called local) quiver data $Q_\xi$, $d_\xi$, $n_\xi$ as follows: the quiver $Q_\xi$ has vertices $1,\ldots,s$ with $\delta_{i,j}-\langle d^i,d^j\rangle$ arrows from $i$ to $j$ for $i,j=1,\ldots,s$. The dimension vector $d_\xi$ is defined by $(d_\xi)_i=m_i$ for $i=1,\ldots,s$, and the dimension vector $n_\xi$ is defined by $(n_\xi)_i=n\cdot d^i$. With this notation, we have the following result (see \cite{SM}):
\begin{theorem}\label{strat} The variety $M_{d,n}^\Theta(Q)$ admits a stratification (in the sense defined above) by the locally closed subsets $M_{d,n}^\Theta(Q)_\xi=\pi_d^{-1}S_\xi$. Each $M_{d,n}^\Theta(Q)_\xi$ admits a fibration (that is, an \'etale locally trivial surjection) over the corresponding Luna stratum $S_\xi$, whose fibre is isomorphic to ${\rm Hilb}_{d_\xi,n_\xi}^{nilp}(Q_\xi)$.
\end{theorem}
By a cell decomposition of a variety $X$ we mean a filtration $\emptyset=X_0\subset X_1\subset\ldots\subset X_s=X$ by closed subvarieties, such that the complements $X_s\setminus X_{s-1}$ are isomorphic to affine spaces.\\[1ex]
For every vertex $i\in I$, we construct a (locally finite) tree quiver $Q_i$ as follows: the vertices $\omega$ of $Q_i$ are indexed by the paths in $Q$ starting in $i$ (including the empty path from $i$ to $i$ of length $0$); there is an arrow $\omega\rightarrow\alpha\omega$ for every path $\omega$ from $i$ to $j$ and every arrow $\alpha:j\rightarrow k$. Note that $Q_i$ has a unique source, corresponding to the empty path. By a subtree $T$ of $Q_i$ we mean a full subquiver which is closed under taking predecessors. The dimension vector $\dimv T$ is defined by setting $(\dimv T)_j$ as the number of paths $\omega\in T$ which end in $j$. By an $n$-forest we mean a tuple $T_*=(T_{i,k})_{i\in I,\, k=1,\ldots,n_i}$ of subtrees $T_{i,k}$ of $Q_i$; its dimension vector is defined as $\dimv T_*=\sum_{i\in I}\sum_{k=1}^{n_i}\dimv T_{i,k}$. It is proved in \cite{SM} that
\begin{theorem}\label{celldec} For all $d$ and $n$, the Hilbert scheme ${\rm Hilb}_{d,n}(Q)$ admits a cell decomposition, whose cells are parametrized by the $n$-forests of dimension vector $d$.
\end{theorem}

\section{Functional equation for $\chi({\rm Hilb}_{d,n}(Q))$ and the big local quiver}\label{fe1}

It follows immediately from Theorem \ref{celldec} that the Euler characteristic of the Hilbert scheme ${\rm Hilb}_{d,n}(Q)$ can be computed as the number of $n$-forests of dimension vector $d$. This allows us to characterize the generating function of these Euler characteristics by a functional equation. For all $n\in \Lambda^+$, we write
$$F^n(t)=\sum_{d\in\Lambda^+}\chi({\rm Hilb}_{d,n}(Q))t^d\in{\bf Q}[[\Lambda]].$$

\begin{proposition}\label{grafting} The series $F^n(t)$ are the uniquely determined elements of ${\bf Q}[[\Lambda]]$ satisfying the following functional equations:
\begin{enumerate}
\item For all $n\in\Lambda^+$, we have $F^n(t)=\prod_{i\in I}F^i(t)^{n_i}$,
\item for all $i\in I$, we have $F^i(t)=1+t_i\prod_{j\in I}F^j(t)^{r_{i,j}}.$
\end{enumerate}
\end{proposition}

\proof Comparing coefficients of $t^d$ in both sides of the first identity, we see that the first claim reduces to the definition of $n$-forests. With the same method, the second identity reduces to the existence of a bijection between subtrees of $Q_i$ of dimension vector $d$ and tuples $(T_{j,k})_{j\in I,\, ,k=1,\ldots,r_{i,j}}$ of subtrees $T_{j,k}$ of $Q_j$ such that $\sum_{j\in I}\sum_{k=1}^{r_{i,j}}\dimv T_{j,k}=d-i$. Such a bijection is provided, by definition of the trees $Q_i$, by grafting the subtrees $T_{j,k}$ to a common root $i$ to obtain any subtree of $Q_i$ exactly once. \hb

\remark In the special case of a quiver with a single vertex and a number of loops, this functional equation is derived in \cite{RNCHilb}.

\begin{proposition} For all $d,n\in\Lambda^+$, we have $\chi({\rm Hilb}_{d,n}^{nilp}(Q))=\chi({\rm Hilb}_{d,n}(Q))$.
\end{proposition}

\proof We adopt an argument used in \cite{CBVdB}. There is a natural ${\bf C}^*$-action on representations of $Q$ by rescaling the maps representing the arrows by a common factor. This action induces actions on ${\rm Hilb}_{d,n}(Q)$ and $M_d^{ssimp}(Q)$, for which the map $\pi_d: {\rm Hilb}_{d,n}(Q)\rightarrow M_d^{simp}(Q)$ is equivariant. Moreover, there exists a unique fixed point for the action of ${\bf C}^*$ on $M_d^{ssimp}(Q)$, namely the point $0$, to which all points of $M_d^{ssimp}(Q)$ attract, in the sense that $\lim_{t\rightarrow 0}t\cdot M=0$ for all $M\in M_d^{ssimp}(Q)$. Therefore, all points of ${\rm Hilb}_{d,n}(Q)$ admit a well-defined limit in the projective variety $\pi^{-1}(0)={\rm Hilb}_{d,n}^{nilp}(Q)$. For each connected component $C$ of ${\rm Hilb}_{d,n}^{nilp}(Q)$, we have its attractor $A_C$ consisting of all points of ${\rm Hilb}_{d,n}(Q)$ whose limit belongs to $C$. By the Bialynicki-Birula theorem \cite{BB}, the attractors $A_C$ are affine fibrations over the components $C$. Consequently, the Euler characteristics of ${\rm Hilb}_{d,n}(Q)$ and of ${\rm Hilb}_{d,n}^{nilp}(Q)$ coincide.\hb

Now we fix data $Q,\Theta,\mu,n$ as before, and associate to it a locally finite quiver (called the big local quiver) $\widetilde{Q}$ as follows: the vertices of $\widetilde{Q}$ are indexed by pairs $(d,i)$ in $\PrL^+_\mu\times{\bf N}$.  The number of arrows from vertex $(d,i)$ to $(d',i')$ is given as $\delta_{d,d'}\cdot\delta_{i,i'}-\langle d,d'\rangle$. For a function $l:\PrL^+_\mu\rightarrow{\bf N}$, we define $\widetilde{Q}_l$ as the full subquiver of $\widetilde{Q}$ supported on the set of vertices $(d,i)$ for $d\in\PrL^+_\mu$ and $1\leq i\leq l(d)$.\\[1ex]
We define dimension vectors $\widetilde{n}$ for the various quivers $\widetilde{Q}_l$ by $\widetilde{n}_{(d,i)}=n\cdot d$. The product $S(\widetilde{Q})=\prod_{d\in\PrL^+_\mu}S_\infty$ of infinite symmetric groups acts on the vertices of $\widetilde{Q}$ by permutation $(\sigma_e)_{e\in\Lambda^+_\mu}(d,i)=(d,\sigma_d(i))$; this restricts to an action of $\prod_{d\in\PrL^+_\mu}S_{l(d)}$ on $\widetilde{Q}_l$.\\[2ex]
For a polystable type $\xi=((d^1,\ldots,d^s),(m_1,\ldots,m_s))$ as above, we can view the local quiver $Q_\xi$ as the quiver $\widetilde{Q}_{l_\xi}$ just defined, where the function $l_\xi$ is given by defining $l_\xi(d)$ as the number of indices $1\leq j\leq s$ such that $d=d^j$. The dimension vector $d_\xi$ for $Q_\xi$ can then be viewed as a dimension vector $\widetilde{d}_\xi$ for $Q_l$. This dimension vector can be made unique by assuming that its entries $(\widetilde{d}_xi)_{d,i}$, for fixed $d\in\PrL_\mu^+$, form a partition, that is, $(\widetilde{d}_\xi)_{(d,1)}\geq\ldots\geq(\widetilde{d}_xi)_{(d,l_\xi(d))}$. Therefore, we call dimension vectors $\widetilde{d}$ of $\widetilde{Q}_l$ partitive if $\widetilde{d}_{(d,1)}\geq\ldots\geq\widetilde{d}_{(d,l(d))}$ for all $d\in\Lambda^+_\mu$; the set of all partitive dimension vectors for $\widetilde{Q}$ (resp. $\widetilde{Q}_l$) is denoted by $\Lambda(\widetilde{Q})^\geq$ (resp. $\Lambda(\widetilde{Q}_l)^\geq$).
We have a natural specialization map $\nu:\Lambda(\widetilde{Q}_l)^+\rightarrow\Lambda^+_\mu$ given by $\nu(d,i)=d$.\\[2ex]
We consider the generating function
$$R^n_l(t)=\sum_{\widetilde{d}\in\Lambda^+(\widetilde{Q}_l)}\chi({\rm Hilb}_{\widetilde{d},\widetilde{n}}(\widetilde{Q}_l)t^{\nu(\widetilde{d})}\in{\bf Z}[[\Lambda^+_\mu]],$$
the specialization of the generating function $F^{\widetilde{n}}$ for the quiver $\widetilde{Q}_l$ with respect to the map $\nu$. By the natural $\prod_{d\in\PrL^+_\mu}S_{l(d)}$-symmetry of $\widetilde{Q}_l$, we have $R^{(d,i)}_l(t)=R^{(d,j)}_l(t)$ for all $d\in\Lambda^+_\mu\setminus 0$ and all $1\leq i,j\leq l(d)$. We denote this series by $R^{(d)}_l(t)$.
Applying Proposition \ref{grafting} and the definition of $\widetilde{Q}_l$, we get
$$R^n_l(t)=\prod_{d\in\PrL^+_\mu}R^{(d)}_l(t)^{l(d)\cdot(n\cdot d)}$$
and
$$R_l^{(d)}(t)=1+t^d\cdot R_l^{(d)}(t)\cdot\prod_{e\in\PrL^+_\mu}R^{(e)}_l(t)^{-\langle d,e\rangle\cdot l(e)}.$$
Call a dimension vector for $\widetilde{Q}_l$ faithful if all its entries are non-zero, and denote by $\Lambda(\widetilde{Q}_l)^{++}$ the set of all such dimension vectors. Define
$$\PrR^n_l(t)=\sum_{\widetilde{d}\in\Lambda^+(\widetilde{Q}_l)^{++}}\chi( {\rm Hilb}_{\widetilde{d},\widetilde{n}}(\widetilde{Q}_l)t^d\in{\bf Z}[[\Lambda^+_\mu]].$$
Using again the symmetry of $\widetilde{Q}_l$, we see that
$$R^n_l(t)=\sum_{l':\PrL^+_\mu\rightarrow{\bf N}}\prod_{d\in\PrL^+_\mu}\binom{l(d)}{l'(d)}\cdot\PrR^n_{l'}(t).$$

Let $\chi:\Lambda^+_\mu\setminus 0\rightarrow{\bf Z}$ be a function with arbitrary integer values (in contrast to the function $l$ considered so far), and define a formal series by
$$R^n_\chi(t)=\sum_{l':\PrL^+_\mu\rightarrow{\bf N}}\prod_{d\in\PrL^+_\mu}{\binom{\chi(d)}{l'(d)}}\cdot\PrR^n_{l'}(t).$$

Simililarly to the above, we have series $\PrR_l^{(d)}(t)$ and $R_\chi^{(d)}(t)$ for $d\in\Lambda^+_\mu$ as special cases of the series $\PrR_l^n(t)$ and $R_\chi^n(t)$, respectively.

\begin{lemma}\label{trick} The series $R^n_\chi(t)$ are given by the functional equations
$$R^n_\chi(t)=\prod_{d\in\PrL^+_\mu}R^{(d)}_\chi(t)^{\chi(d)\cdot(n\cdot d)}$$
and
$$R_\chi^{(d)}(t)=1+t^d\cdot R_\chi^{(d)}\cdot\prod_{e\in\PrL^+_\mu}R^{(e)}_\chi(t)^{-\langle d,e\rangle\cdot \chi(e)}.$$
\end{lemma}

\proof It is easy to see that there exist unique series $S_\chi^d(t)$ for all functions $\chi$ as above and all $d\in\PrL^+_\mu$ fulfilling the equations
$$S_\chi^d(t)=1+t^d\cdot S_\chi^d(t)\cdot\prod_{e\in\PrL^+_\mu}S^e_\chi(t)^{-\langle d,e\rangle\cdot \chi(e)},$$
since these functional equations define recursions determining the coefficients of the series.
These coefficients depend polynomially on the values $\chi(d)$. The same holds for the coefficients of the series $R_\chi^{(d)}(t)$ by definition. Now the equality $S_\chi^d(t)=R_\chi^{(d)}(t)$ holds for all functions $\chi$ with values in ${\bf N}$, thus it has to hold for arbitrary $\chi$.\hb

\section{Functional equation for $\chi(M_{d,n}^\Theta(Q))$}\label{section4}

We start with a calculation of Euler characteristics of strata of symmetric products of a variety, which should be well-known.
Denote by $\mathcal{P}$ the set of all partitions. For $\lambda$ in $\mathcal{P}$, denote by $m_i(\lambda)$ the multiplicity of $i$ in $\lambda$, that is, the number of indices $j$ such that $\lambda_j=i$. For a variety $X$, we denote by $S^nX$ its $n$-th symmetric power, that is, the quotient of $X^n$ by the natural action of the symmetric group $S_n$. The product variety $X^n$ admits a stratification by strata $X^n_I$, where $I=(I_1,\ldots,I_k)$ is a decomposition of $\{1,\ldots,n\}$ into pairwise disjoint subsets. Namely, $X^n_I$ is defined as the set of ordered tuples $(x_1,\ldots,x_n)$ such that $x_i=x_j$ if and only if $i,j$ belong to the same subset $I_l$. Obviously, $X^n_I$ is isomorphic to $X^k_{(1,\ldots,1)}$, the set of unordered $k$-tuples of pairwise different points in $X$.

Any $I$ as above induces a partition $\lambda(I)$ of $n$, with parts being the cardinalities of the subsets $I_k$ forming $I$. The image of $X^n_I$ under the quotient map $\pi:X^n\rightarrow S^nX$ depends only on the partition $\lambda=\lambda(I)$ and is denoted by $S^n_\lambda X$. The inverse image under $\pi$ of $S^n_\lambda X$ is precisely the union of the strata $X^n_I$ such that $\lambda(I)=\lambda$. Moreover, the fibre of $\pi$ over a point in $S^n_\lambda X$ is finite of cardinality $\frac{n!}{\lambda_1!\ldots\lambda_k!}$. The number of decompositions $I$ such that $\lambda(I)=\lambda$ equals
$$\frac{n!}{\lambda_1!\cdot\ldots\cdot\lambda_k!}\cdot\frac{1}{\prod_i(m_i(\lambda)!)}.$$
An easy induction shows that the Euler characteristic in cohomology with compact support $\chi$ of $X^n_{(1,\ldots,1)}$ equals $$\chi(X)(\chi(X)-1)\ldots(\chi(X)-n+1)={n!}{\binom{\chi(X)}{n}}.$$
We have thus proved:

\begin{lemma}\label{lemmasymm} For all partitions $\lambda$ of $n$, we have $$\chi(S^n_\lambda X)=\frac{1}{\prod_im_i(\lambda)!}\chi(X)(\chi(X)-1)\ldots(\chi(X)-k+1)=\frac{1}{\prod_im_i(\lambda)!}k!{\binom{\chi(X)}{k}}.$$
\end{lemma}

We can now consider the generating function of the Euler characteristics of arbitrary smooth models, using the big local quiver notation of the previous section.\\[2ex]
In particular, to a polystable type $\xi$, we have associated a partitive dimension vector $p$ for $\widetilde{Q}$ (resp. a large enough $\widetilde{Q}_l$); we denote the stratum $S_\xi$ by $S_p$. With the above notation, we have $$S_p\simeq\prod_{d\in\PrL^+_\mu}S_{p(d)}^{|p(d)|}M_d^{st}(Q)$$
by definition of $S_\xi$. Theorem \ref{strat} can now be rephrased as stating that $M_{d,n}^\Theta(Q)$ admits a stratification indexed by partitive dimension vectors $p\in\Lambda(\widetilde{Q})^+$ such that $\nu(p)=d$. Each stratum is a locally trivial fibration over $S_p$, with fibre isomorphic to ${\rm Hilb}_{p,\widetilde{n}}^{nilp}(\widetilde{Q})$.
We thus have, using Lemma \ref{lemmasymm} for the second equality:
\begin{eqnarray*}
\chi(M_{d,n}^\Theta(Q))&=&\sum_p\chi(S_p)\cdot\chi({\rm Hilb}_{p,\widetilde{n}}^{nilp}(\widetilde{Q}))\\
&=&\sum_p\prod_{d\in\PrL_\mu^+}\frac{1}{\prod_im_i(p(d))!}l(p(d))!{\binom{\chi(M_d^{st}(Q))}{l(p(d))}}\cdot\chi({\rm Hilb}_{p,\widetilde{n}}^{nilp}(\widetilde{Q})),
\end{eqnarray*} 
the sum running over all partitive dimension vectors $p$ for $\widetilde{Q}$ such that $\nu(p)=d$.\\[2ex]
Considering the generating function, we thus have
$$\sum_{d\in\Lambda^+_\mu}\chi(M_{d,n}^\Theta(Q))t^d$$
$$=\sum_{p\in\Lambda(\widetilde{Q})^\geq}\prod_{d\in\PrL^+_\mu}\left(\frac{1}{\prod_im_i(p(d))!}l(p(d))!{\binom{\chi(M_d^{st}(Q))}{l(p(d))}}\right)\cdot\chi({\rm Hilb}_{p,\widetilde{n}}^{nilp}(\widetilde{Q}))t^{\nu(p)}.$$
Sorting by lengths of the partitions, this can be rewritten as
$$\sum_{l:\PrL^+_\mu\rightarrow{\bf N}}\sum_{p\in\Lambda(\widetilde{Q}_l)^\geq}\prod_{d\in\PrL^+_\mu}\left(\frac{1}{\prod_im_i(p(d))!}l(d))!{\binom{\chi(M_d^{st}(Q))}{l(d))}}\right)\cdot\chi({\rm Hilb}_{p,\widetilde{n}}^{nilp}(\widetilde{Q}))t^{\nu(p)}.$$
We want to extend the range of summation in the inner sum to arbitrary dimension vectors  for each $\widetilde{Q}_l$ without changing the sum. By the symmetry property of $\widetilde{Q}$ (resp.~$\widetilde{Q}_l$), we can do this by incorporating a factor which counts the number of derangements of a given partitive dimension vector $p$ into arbitrary dimension vectors. This number is precisely
$$\prod_{d\in\PrL^+_\mu}\frac{l(p(d))!}{\prod_im_i(p(d))!},$$
this factor being already present. Thus, the above sum equals
$$\sum_{l:\PrL^+_\mu\rightarrow{\bf N}}\sum_{\widetilde{d}\in\Lambda(\widetilde{Q}_l)^{++}}\prod_{d\in\PrL^+_\mu}{\binom{\chi(M_d^{st}(Q)}{l(d)}}\cdot\chi({\rm Hilb}_{\widetilde{d},\widetilde{n}}^{nilp}(\widetilde{Q}))t^{\nu(\widetilde{d})},$$
the inner sum now running over all faithful dimension vectors for $\widetilde{Q}_l$. Using the previous notation, this equals
$$\sum_{l:\PrL^+_\mu\rightarrow {\bf N}}\prod_{d\in\PrL_\mu^+}{\binom{\chi(M_d^{st}(Q))}{l(d)}}\PrR_l^n(t)=R^n_\chi(t)$$
for the function $\chi$ defined by $\chi(d)=\chi(M_d^{st}(Q))$. By Lemma \ref{trick}, we arrive at the following result:

\begin{theorem}\label{t42} The generating function of Euler characteristics of smooth models is defined by the functional equations
$$\sum_{d\in\Lambda^+_\mu}\chi(M_{d,n}^\Theta(Q))t^d=\prod_{d\in\PrL^+_\mu}R^{d}(t)^{\chi(M_d^{st}(Q))\cdot(n\cdot d)}$$
and
$$R^{d}(t)=1+t^d\cdot R^d(t)\cdot\prod_{e\in\PrL^+_\mu}R^{e}(t)^{-\langle d,e\rangle\cdot \chi(M_e^{st}(Q))}.$$
\end{theorem}

To make the nature of these functional equations more transparent, we will define a slight variant of the generating functions. Writing $$Q^n_\mu(t)=\sum_{d\in\Lambda_\mu^+}\chi(M_{d,n}^\Theta(Q))t^d,$$ we have $Q^n_\mu(t)=\prod_{i\in I}Q^{i}_\mu(t)^{n_i}$
by the previous theorem. This suggests the definition $Q^\eta_\mu(t)=\prod_{i\in I}Q^{i}_\mu(t)^{\eta(i)}$
for an arbitrary linear functional $\eta\in\Lambda^*$, so that $Q^{n\cdot}_\mu(t)=Q^n_\mu(t)$ for all $n\in\Lambda^+$. In particular, we consider $S^d_\mu(t)=Q^{\langle d,\_\rangle}_\mu(t)$ for $d\in\PrL^+_\mu$.

\begin{corollary}\label{corsd} The series $S^d_\mu(t)$ for $d\in\PrL^+_\mu$ are given by the functional equations
$$S^d_\mu(t)=\prod_{e\in\PrL_\mu^+}(1-t^eS^e_\mu(t))^{-\langle d,e\rangle\cdot\chi(M_e^{st}(Q))}.$$
\end{corollary}

\proof By the definitions and Theorem \ref{t42}, we have
$$S^d_\mu(t)=\prod_{e\in\PrL_\mu^+}R^e(t)^{\langle d,e\rangle\cdot\chi(M_e^{st}(Q))}.$$
The last line of Theorem \ref{t42} can be restated as
$$R^d(t)=(1-t^d\prod_{e\in\PrL_\mu^+}R^e(t)^{-\langle d,e\rangle\cdot\chi(M_e^{st}(Q))})^{-1},$$
thus
$$R^d(t)=(1-t^dS^d_\mu(t))^{-1}.$$
Substituting this in the factorization of $S^d_\mu(t)$ yields the desired equation.\hb

\section{Duality for Euler products}\label{number}

Let $F(t)\in{\bf Q}[[t]]$ be a formal power series with constant term $F(0)=1$. Then we can write $F(t)$ as an Euler product \begin{equation}\label{eulerproduct} F(t)=\prod_{i\geq 1}(1-(-t)^i)^{-ia_i}\end{equation} for $a_i\in{\bf Q}$ (note the sign convention, which is essential in the following; see the example at the end of this section). We can also characterize $F(t)$ as the unique solution of a functional equation of the form
\begin{equation}\label{functionalequation} F(t)=\prod_{i\geq 1}(1-(tF(t))^i)^{ib_i}\end{equation}
for $b_i\in{\bf Q}$; see the remark below for the proof.\\[1ex]
The main result of this section is:

\begin{theorem}\label{duality} In the above notation, we have $b_i\in{\bf Z}$ for all $i\geq 1$ if and only if $a_i\in{\bf Z}$ for all $i\geq 1$.
\end{theorem}

\remark Writing $H(t)=-tF(t)$, we have, by a straightforward calculation,
$$H(t)=-t\prod_{i\geq 1}(1-(-t)^i)^{-ia_i}$$
and
$$t=-H(t)\prod_{i\geq 1}(1-(-H(t))^i)^{-ib_i}.$$
This means that $H(t)$ is the compositional inverse of the series $$-t\prod_{i\geq 1}(1-(-t)^i)^{-ib_i}.$$
This shows that the series $F(t)$ can be characterized by a functional equation of the form (\ref{functionalequation}) for unique $b_i$, and it shows the symmetry of the statement in the theorem. Thus, we only have to prove integrality of the $a_i$ given integrality of the $b_i$.\\[1ex]
As the first step towards the proof of the theorem, we will derive an explicit formula for the $a_i$ in terms of the $b_i$ by applying Lagrange inversion to the functional equation (\ref{functionalequation}. We use the following version of Lagrange inversion:
\begin{lemma}\label{lagrangeinversion}
Suppose that power series $F(t)$, $G(t)\in{\bf Q}[[t]]$ with $G(0)\not=0$ are related by $F(t)=G(tF(t))$. Then, for all $k,d\in{\bf Z}$, we have
$$(k+d)[t^d]F(t)^k=k[t^d]G(t)^{k+d},$$
where $[t^d]F(t)$ denotes the $t^d$-coefficient of the series $F(t)$.
\end{lemma}

\proof Apply \cite[Theorem 5.4.2]{St} using the notation $f(t)=tF(t)$ and $d=n-k$.\hb

\begin{lemma}\label{coeffpartition} For all $d\in{\bf N}$ and all $c_i\in{\bf Z}$ for $i\geq 1$, we have
\begin{equation}[t^d]\prod_{i\geq 1}(1-t^i)^{-c_i}=\sum_{\lambda\vdash d}\prod_{i\geq 1}\binom{c_i+\lambda_i-\lambda_{i+1}-1}{\lambda_i-\lambda_{i+1}},\end{equation}
the sum ranging over all partitions $\lambda$ of $d$.
\end{lemma}

\proof We have $$(1-t)^{-c}=\sum_{k\geq 0}\binom{c+k-1}{k}t^k,$$
and therefore
\begin{eqnarray*}[t^d]\prod_{i\geq 1}(1-t^i)^{-c_i}&=&[t^d]\prod_{i\geq 1}\sum_{k_i\geq 0}\binom{c_i+k_i-1}{k_i}t^k_i=\\
&&[t^d]\sum_{k_1,k_2,\ldots\geq 0}\prod_{i\geq 1}\binom{c_i+k_i-1}{k_i}t^{\sum_ik_i}=\\
&&[t^d]\sum_{\lambda}\prod_{i\geq 1}\binom{c_i+\lambda_i-\lambda_{i+1}-1}{\lambda_i-\lambda_{i+1}}t^{|\lambda},\end{eqnarray*}
where the last sum ranges over all partitions $\lambda$, which are related to sequences $k_1,k_2,\ldots\geq 0$ via $\lambda_i=\sum_{j\geq i}k_j$.\hb

\remark Here and in the following, we make frequent use of binomial coefficients $\binom{a}{b}$ for $a\in{\bf Z}$ using
\begin{equation}\label{negativebinom}\binom{-a+b-1}{b}=(-1)^b\binom{a}{b}\end{equation}

Using these preparations, we can state the desired formula relating the coefficients $a_i$ and $b_i$:

\begin{proposition}\label{moebiusinversion} With the above notation, we have, for all $d\geq 1$:
\begin{equation}\label{mif}d^2a_d=\sum_{e|d}\mu(d/e)(-1)^e\sum_{\lambda\vdash e}(-1)^{\lambda_1}\prod_{i\geq 1}\binom{ib_ie}{\lambda_i-\lambda_{i+1}},\end{equation}
where the first sum ranges over all divisors of $d$, and $\mu$ denotes the number-theoretic Moebius function.
\end{proposition}

\proof We apply Lemma \ref{lagrangeinversion} to the functional equation (\ref{functionalequation}) using
$$G(t)=\prod_{i\geq 1}(1-t^i)^{ib_i}$$ and get
\begin{equation}\label{lifn}(k+d)[t^d]\prod_{i\geq 1}(1-(-t)^i)^{-ia_ik}=k[t^d]\prod_{i\geq 1}(1-t^i)^{ib_i(k+d)}.\end{equation}
Lemma \ref{coeffpartition} allows us to write the left hand side of (\ref{lifn}) as
$$(k+d)(-1)^d\sum_{\lambda\vdash d}\prod_{i\geq 1}\binom{ia_ik+\lambda_i-\lambda_{i+1}-1}{\lambda_i-\lambda_{i+1}},$$
and the right hand side of (\ref{lifn}) as
$$k\sum_{\lambda\vdash d}\prod_{i\geq 1}\binom{-ib_i(k+d)+\lambda_i-\lambda_{i+1}-1}{\lambda_i-\lambda_{i+1}}.$$
We use (\ref{negativebinom}) and substitute $k$ by $X$ to rewrite (\ref{lifn}) as
\begin{equation}\label{rewrite}X\sum_{\lambda\vdash d}(-1)^{\lambda_1}\prod_{i\geq 1}\binom{ib_i(X+d)}{\lambda_i-\lambda_{i+1}}=(-1)^d(X+d)\sum_{\lambda\vdash d}\prod_{i\geq 1}\binom{ia_iX+\lambda_i-\lambda_{i+1}-1}{\lambda_i-\lambda_{i+1}}.\end{equation}
Both sides behaving polynomially in $X$, equality for all $X\in{\bf Z}$ thus implies equality of the polynomials. We want to compare the linear $X$-terms (the constant terms being $0$) of both sides. Note the following property:\\[1ex]
The polynomial $\binom{aX+b+c-1}{c}$
has constant $X$-coefficient $\binom{b+c-1}{c}$, and the polynomial $\binom{aX+c-1}{c}$ has linear $X$-coefficient $a/c$.\\[1ex]
Applying this, we see that the left hand side of (\ref{rewrite}) has linear $X$-coefficient
$$\sum_{\lambda\vdash d}(-1)^{\lambda_1}\prod_{i\geq 1}\binom{ib_id}{\lambda_i-\lambda_{i+1}}.$$
To analyze the linear $X$-coefficient of the right hand side of (\ref{rewrite}), note first that the constant $X$-coefficient of each product
\begin{equation}\label{product}\prod_{i\geq 1}\binom{ia_iX+\lambda_i-\lambda_{i+1}-1}{\lambda_i-\lambda_{i+1}}\end{equation}
equals zero. Its linear $X$-term is non-zero only if exactly one factor appears, that is, if there is only one non-zero difference $\lambda_i-\lambda_{i+1}$. In this case, the partition $\lambda$ of $d$ equals $$\lambda=\underbrace{(d/i,\ldots,d/i)}_{i\mbox{-times}}$$
for a divisor $i$ of $d$. Thus, the  product (\ref{product}) reduces to
$$\binom{ia_iX+d/i-1}{d/i},$$
having linear $X$-coefficient $(ia_i)/(d/i)=i^2a_i/d$ by the above. We conclude that the linear $X$-coefficient of the right hand side of (\ref{rewrite}) equals
$$(-1)^d\sum_{i|d}i^2a_i.$$
Comparison of both linear $X$-coefficients thus yields
$$\sum_{i|d}i^2a_i=(-1)^d\sum_{\lambda\vdash d}(-1)^{\lambda_1}\prod_{i\geq 1}\binom{ib_id}{\lambda_i-\lambda_{i+1}}.$$
After Moebius inversion, we arrive at the claimed formula (\ref{mif}).\hb

To prove integrality of the $a_d$ given integrality of all $b_i$, we thus have to prove that the right hand side of (\ref{mif}) is divisible by $d^2$. This can be tested on the prime divisors of $d$. Denoting by $$m(d)=m_p(d)=\max\{m\,:\, p^m|d\}$$ the multiplicity of a prime $p$ as a divisor of $d$, we thus have to prove divisibility by $p^{2m_p(d)}$ of the right hand side of (\ref{mif}) for all primes $p$. We prepare this proof by stating certain divisibility/congruence properties of binomial coefficients.

\begin{lemma}\label{kummer} Let $p$ be a prime. For $a,b\in{\bf Z}$ and $b\geq 0$, we have
$$p^{\max(m_p(a)-m_p(b),0)}|\binom{a}{b}.$$
\end{lemma}

\proof By a result of Kummer (see, for example, \cite{Gr}), the exact power of $p$ diving $\binom{a}{b}$ equals the number of `carries' when subtracting $b$ from $a$ in base $p$, at least when $a\geq 0$. This can be generalized to $a\in{\bf Z}$ using
\begin{equation}\label{negbinom}\binom{-a}{b}=(-1)^b\frac{a}{a+b}\binom{a+b}{b}.\end{equation}
The lemma follows.\hb

\begin{lemma}\label{gessel} Let $p$ be a prime, and define $\mu_p=0,1,2$ provided $p=2$, $p=3$, $p\geq 5$, respectively. Assume $p|a,b$ for integers $a$, $b$ with $b\geq 0$. Define $\eta$ as $-1$ if $p=2$ and $b\equiv 2\equiv a-b\bmod 4$, and as $1$ otherwise. Then
$$\binom{a}{b}\equiv\eta\binom{a/p}{b/p}\bmod p^r,$$
for
$$r\leq m_p(a)+m_p(b)+m_p(a-b)+m_p(\binom{a/p}{b/p})-\mu_p.$$
In case $p=2$, we also have $$\binom{a}{b}\equiv\binom{a/2}{b/2}\bmod 4.$$
\end{lemma}

\proof The general statement (usually \cite{Ge,Gr} attributed to Jacobsthal \cite{Jac}) is proved in \cite[Theorem 2.2]{Ge}, with the assumption $a\geq 0$ there removed by (\ref{negbinom}). For the congruence modulo $4$, we calculate as in the proof of \cite[Theorem 2.2]{Ge}:
$$\binom{a}{b}=\binom{a/2}{b/2}\prod_{\stackrel{i=1}{2\nmid i}}^b(1+2(a-b)/i)\equiv\binom{a/2}{b/2}(1+2(a-b)\sum_{\stackrel{i=1}{2\nmid i}}^b 1/i)\equiv$$
$$\equiv\binom{a/2}{b/2}(1+(a-b)(b/2)^2)\bmod 4.$$
The term $(a-b)(b/2)^2$ is congruent to $1 \bmod 4$ except when $b/2$ is odd and $a/2$ is even, in which case it is congruent to $-1 \mod 4$. But in this case, $\binom{a/2}{b/2}$ is even by Lemma \ref{kummer}.\hb

From the previous two lemmas, we derive divisibility/congruence properties of the product of binomial coefficients appearing in (\ref{mif}).

\begin{lemma}\label{claim2} Let $p$ be a prime dividing $e\geq 0$. If $\lambda$ is a partition of $e$ which is not divisible by $p$ (that is, some coefficient $\lambda_i$ is not divisible by $p$), we have
$$p^{2m(e)}|\prod_{i\geq 1}\binom{ib_ie}{\lambda_i-\lambda_{i+1}}.$$
\end{lemma}

\proof To shorten notation, we write $m=m_p(e)$ and $c_i=\lambda_i-\lambda_{i+1}$ for $i\geq 1$, thus $e=\sum_{i\geq 1}ic_i$. Lemma \ref{kummer} yields
$$p^{\max(m+m(i)-m(c_i),0)}|\binom{ib_ie}{c_i};$$
we thus have to prove
\begin{equation}\label{ineq}\sum_{i:c_i\not=0}\max(m+m(i)-m(c_i),0)\geq 2m\end{equation}
provided some $c_i\not=0$ is not divisible by $p$. Let $i_0$ be an index such that $m(c_{i_0})=0$.\\[1ex]
Let $m_0$ be the minimum over all $m(i)+m(c_i)$. Since $e=\sum_iic_i$, we can distinguish two cases: either $m_0=m$ (case 1), or $m_0<m$ and the minimum is obtained at least twice (case 2). For case 1 we have, in particular, $m(i_0)\geq m$, thus $$\max(m+m(i_0)-m(c_{i_0}),0)\geq 2m,$$
and (\ref{ineq}) follows.\\[1ex]
For case 2, let $i_1,i_2$ be two different indices where the minimum $m_0$ is obtained. For $s=1,2$, we have $m+m(i_s)-m(c_{i_s})\geq 0$, since otherwise,
$$m>m_0=m(c_{i_s})+m(i_s)\geq m(c_{i_s})>m+m(i_s),$$
a contradiction.  Again we distinguish two cases: first, assume that $i_0$ coincides with, say, $i_1$. Then we can estimate \begin{eqnarray*}&&\sum_{i:c_i\not=0}\max(m+m(i)-m(c_i),0)\\
&\geq&\max(m+m(i_0)-m(c_{i_0}),0)+\max(m+m(i_2)-m(c_{i_2}),0)\\
&=&2m+m_0+m(i_2)-m(c_{i_2})=2m+2m(i_2)\geq 2m,\end{eqnarray*}
and (\ref{ineq}) follows. Second, assume that $i_0$ differs from $i_1$, $i_2$. Like in the previous case, we can estimate \begin{eqnarray*}&&\sum_{i:c_i\not=0}\max(m+m(i)-m(c_i),0)\\
&\geq&3m+m(i_0)+m(i_1)+m(i_2)-m(c_{i_1})-m(c_{i_2})\\
&\geq&2m+m-m_0+2m(i_1)+2m(i_2)\geq 2m,\end{eqnarray*}
and (\ref{ineq}) follows again.\hb

\begin{lemma}\label{claim1} Let $p$ be a prime dividing $e\geq 0$. If $\lambda=p\mu$ is a partition of $e$ divisible by $p$, then
\begin{equation}\label{cong}\prod_{i\geq 1}\binom{ib_ie}{\lambda_i-\lambda_{i+1}}\equiv (-1)^{(p-1)(e/p+\mu_1)}\prod_{i\geq 1}\binom{ib_ie/p}{\mu_i-\mu_{i+1}}\bmod p^{2m_p(e)}.\end{equation}
\end{lemma}

\proof So assume that $\lambda=p\mu$, and denote again $m=m(e)$ and $c_i=\lambda_i-\lambda_{i+1}$. Applying the general congruence of Lemma \ref{gessel} to a non-trivial (that is, $c_i\not=0$) factor of the left hand side of (\ref{cong}), we get
$$\binom{ib_ie}{c_i}\equiv\eta_i\binom{ib_ie/p}{c_i/p}\bmod p^{r_i},$$
where the sign $\eta_i$ is $-1$ only in case $p=2$, $c_i/2$ odd, $ib_ie/2-c_i/2$ odd, and
\begin{eqnarray}\label{est}\nonumber r_i&=&m(ib_ie)+m(c_i)+m(ib_ie-c_i)+m(\binom{ib_ie/p}{c_i/p}-\mu_p\\
\nonumber &\geq&m(e)+m(i)+m(c_i)+\min(m(e)+m(i),m(c_i))+\\
\nonumber &&\max(m(e)+m(i)-m(c_i),0)-\mu_p\\
&=&2m+2m(i)+m(c_i)-\mu_p.\end{eqnarray}
Suppose first that $p\geq 3$. Then $r_i\geq 2m$ using $m(c_i)\geq 1$ and $\mu_p\leq 1$. The sign in (\ref{cong}) vanishes due to the even factor $p-1$, and $\eta_i=1$. The congruence (\ref{cong}) follows.\\[1ex]
Next, assume that $p=2$ and $m\geq 2$. Then the estimate (\ref{est}) only assures congruence of the binomial coefficients $\bmod 2^{2m-1}$ in case $i$ and $c_i/2$ are odd, thus
$$\binom{ib_ie}{c_i}\equiv\eta_i\binom{ib_ie/2}{c_i/2}+\varepsilon_i\bmod 2^{2m},$$
where $\varepsilon_i\in\{0,2^{2m-1}\}$, non-triviality only being possible if $i$ and $c_i/2$ are odd. Then
\begin{eqnarray}\label{complex}\nonumber \prod_{i\geq 1}\binom{ib_ie}{c_i}&\equiv&\prod_{i\geq 1}(\eta_i\binom{ib_ie}{c_i/2}+\varepsilon_i)\\
&\equiv&\prod_{i\geq 1}\eta_i\binom{ib_ie/2}{c_i/2}+\sum_{i\geq 1}\varepsilon_i\prod_{j\not=i}\eta_j\binom{jb_je}{c_j/2}\bmod 2^{2m},\end{eqnarray}
since all multiple products of the $\varepsilon_i$ vanish $\bmod 2^{2m}$. For the same reason, we only have to consider summands in (\ref{complex}) for which $\varepsilon_i\not=0$ and each factor $$\binom{jb_je/2}{c_j/2}$$ is odd. Since $m\geq 2$, this can only happen (using Lemma \ref{kummer}) in the case that $m(c_j)\geq m(j)+m$ for all $j\not=i$ such that $c_j\not=0$. But then $$2^m|e-\sum_{j\not=i:c_j\not=0}jc_j=ic_i,$$
a contradiction to the assumptions $m(ic_i)=1$ (by $\varepsilon_i\not=0$) and $m\geq 2$. Thus, we have proved that
$$\prod_{i\geq 1}\binom{ib_ie}{c_i}\equiv\prod_{i\geq 1}\eta_i\cdot\prod_{i\geq 1}\binom{ib_ie/2}{c_i/2}\bmod 2^{2m},$$
and we have to compare the sign $\prod_i\eta_i=(-1)^u$ to the sign of (\ref{cong}). Using $m\geq 2$, we have
\begin{eqnarray*}
u&=&{|\{i\geq 1\, :\, c_i/2\mbox{ odd}, ib_ie/2-c_i/2\mbox{ odd}\}|}\\
&=&{|\{i\geq 1\, :\, c_i/2\mbox{ odd}\}|}.\end{eqnarray*}
The sign in (\ref{cong}) equals
$$(-1)^{e/2+\sum_ic_i/2},$$
and we are done.\\[1ex]
Finally, consider the case $p=2$ and $m=1$. Then the statement on congruences $\bmod 4$ of Lemma \ref{gessel} yields
$$\prod_{i\geq 1}\binom{ib_ie}{c_i}\equiv\prod_{i\geq 1}\binom{ib_ie/2}{c_i/2}\bmod 4,$$
and again we only have to consider the sign.  The sign in (\ref{cong}) equals $$(-1)^{1+\sum_i c_i/2}.$$
We have $e/2=\sum_iic_i/2$, thus the sum $\sum_{{2}{\nmid}\,{i}}c_i/2$ is odd. Suppose $\sum_ic_i/2$ is even (the only case in which the sign of (\ref{cong}) potentially differs from $1$). Then $\sum_{2|i}c_i/2$ is odd. Thus, there exists an even index $i$ with $c_i/2$ odd. In this case, the binomial coefficient $$\binom{ib_ie/2}{c_i/2}$$ is even, and the sign is irrelevant $\bmod 4$.\hb

With these preparations, we can finish the\\[1ex]
{\bf Proof of Theorem \ref{duality}:} Assume that $p$ is a prime such that $m=m(d)=m_p(d)\geq 1$. The divisors $e$ of $d$ for which $\mu(d/e)$ is non-zero fulfill $m(e)=m(d)$ or $m(e)=m(d)-1$, that is, they are of the form $e$ or $e/p$ for a divisor $e$ of $d$ such that $m(e)=m(d)$. We can thus split the right hand side of (\ref{mif}) into the following difference:
\begin{eqnarray}\label{diff}\nonumber &&\sum_{e|d\, :\, m(e)=m(d)}\mu(d/e)(-1)^{e}\sum_{\lambda\vdash e}(-1)^{\lambda_1}\prod_{i\geq 1}\binom{ib_ie}{\lambda_i-\lambda_{i+1}}\\
&-&\sum_{e|d\, :\, m(e)=m(d)}\mu(d/e)(-1)^{e/p}\sum_{\mu\vdash e/p}(-1)^{\mu_1}\prod_{i\geq 1}\binom{ib_ie/p}{\mu_i-\mu_{i+1}}.\end{eqnarray}
Now consider a summand of the first sum of (\ref{diff}) corresponding to a partition $\lambda$ of $e$. If $\lambda$ is not divisible by $p$, then Lemma \ref{claim2} shows that the summand is divisible by $p^{2m(e)}=p^{2m(d)}$. If $\lambda=p\mu$ is divisible by $p$, then Lemma \ref{claim1} shows that the summand is congruent $\bmod p^{2m(d)}$ to the summand of the second sum of (\ref{diff}) corresponding to the partition $\mu$. In other words, the difference of the two sums in (\ref{diff}) vanishes $\bmod p^{2m(d)}$, proving the theorem.\hb

For the application to the integrality of certain Donaldson-Thomas type invariants in the following section, we need a slight generalization of Theorem \ref{duality}. We treat this case separately, although a second inspection of the proofs leading to Theorem \ref{duality} is neccessary, to avoid additional complications in the notation used so far.

\begin{theorem}\label{genduality} Let $F(t)\in{\bf Q}[[t]]$ be a power series with $F(0)=1$. For $N\in{\bf Z}$, write
$$F(t)=\prod_{i\geq 1}(1-((-1)^Nt)^i)^{-ia_i}$$
for $a_i\in{\bf Q}$. We can characterize $F(t)$ as the solution to a functional equation of the form
$$F(t)=\prod_{i\geq 1}(1-(tF(t)^N)^i)^{ib_i}$$
for unique $b_i\in{\bf Q}$. Under these assumptions, we have $b_i\in{\bf Z}$ for all $i\geq 1$ if and only if $a_i\in{\bf Z}$ for all $i\geq 1$.
\end{theorem}

\proof The argument used in the remark following Theorem \ref{duality}, using the power series $H(t)=t(-F(t))^N$, shows existence and uniqueness of the $b_i$, as well as the symmetry of the statement of Theorem \ref{genduality}. Applying Proposition \ref{moebiusinversion} to $G(t)=F(t)^N$ yields the following explicit formula for all $d\geq 1$:
\begin{equation}\label{gmif}d^2a_d=\frac{1}{N}\sum_{e|d}\mu(d/e)(-1)^{Ne}\sum_{\lambda\vdash e}(-1)^{\lambda_1}\prod_{i\geq 1}\binom{Nib_ie}{\lambda_i-\lambda_{i+1}}.\end{equation}
Now any summand of (\ref{gmif}) is divisible by $N$, thus the denominator $N$ in (\ref{gmif}) cancels. Next, note that none of our arguments (Lemma \ref{claim2}, \ref{claim1}) for the proof of Theorem \ref{duality} uses any divisibility properties of the $b_i$, thus these arguments are valid when replacing $b_i$ by $Nb_i$, yielding an additional divisibility by $N$.\\[1ex]
The only additional difficulty is the sign in the statement of Lemma \ref{claim1}, which now reads
$$(-1)^{(p-1)(Ne/p+\mu_1)}.$$
Repeating the sign considerations in the proof of Lemma \ref{claim1}, we see that we can concentrate on the case $p=2$ and $m(e)=1$, where the sign now reads
$$(-1)^{N+\sum_i c_i/2}.$$
The argument of the proof of Lemma \ref{claim1} is still valid in case $N$ is odd. On the other hand, if $N$ is even, we can choose an index $i$ such that $c_i/2$ is odd, and Lemma \ref{kummer} shows that the binomial coefficient
$$\binom{Nib_ie/2}{c_i/2}$$
is even, the sign thus being again irrelevant $\bmod 4$. \hb

\example We consider the example $b_i=0$ for all $i\geq 2$ and denote $b=b_1$. Then $F(t)$ is the solution to the functional equation
$$F(t)=(1-tF(t)^N)^b,$$
and we want to factor $F(t)$ as
$$F(t)=\prod_{i\geq 1}(1-((-1)^Nt)^i)^{-ia_i}.$$
The formula (\ref{gmif}) gives
$$a_d=\frac{1}{Nd^2}\sum_{e|d}\mu(d/e)(-1)^{(N+1)e}\binom{Nbe}{e}.$$
In particular, we have $a_1=(-1)^{N+1}b$ and
$$a_2=\frac{b(2Nb-(1+(-1)^{N+1})}{4},$$
and we see that the choice of signs is essential for the integrality of the $a_d$ given by Theorem \ref{genduality}.\\[1ex]
The particular case $N=1$, $b=-1$ gives a factorization (\ref{eulerproduct}) for the generating function
$$F(t)=\frac{1-\sqrt{1-4t}}{2t}$$
of Catalan numbers with
$$a_d=\frac{1}{d^2}\sum_{e|d}(-1)^e\mu(d/e)\binom{2e-1}{e},$$
which is (up to signs) sequence A131868 in \cite{OEIS}.

\section{Application to Donaldson-Thomas invariants and wall-crossing formulas}\label{app}

In this section, we apply the results of the previous sections to the setup of \cite{RWC}. We assume that $Q$ is a quiver without oriented cycles, thus we can order the vertices as $I=\{i_1,\ldots,i_r\}$ in such a way that $k>l$ provided there exists an arrow $i_k\rightarrow i_l$. Denote by $\{\_,\_\}$ the skew-symmetrization of $\langle\_,\_\rangle$, thus $\{d,e\}=\langle d,e\rangle-\langle e,d\rangle$. Define $b_{ij}=\{i,j\}$ for $i,j\in I$.\\[1ex]
We consider the formal power series ring $B={\bf Q}[[\Lambda^+]]={\bf Q}[[x_i: i\in I]]$ with topological basis $x^d=\prod_{i\in I}x_i^{d_i}$ for $d\in\Lambda^+$. The algebra $B$ becomes a Poisson algebra via the Poisson bracket
$$\{x_i,x_j\}=b_{ij}x_ix_j\mbox{ for }i,j\in I.$$
Define Poisson automorphisms $T_i$ of $B$ by
$$T_i(x_j)=x_j\cdot(1+x_i)^{\{i,j\}}$$
for all $i,j\in I$.\\[1ex]
We study a factorization property in the group ${\rm Aut}(B)$ of Poisson automorphisms of $B$ involving a descending product $\prod^\leftarrow_{\mu\in{\bf Q}}$ indexed by rational numbers, which is indeed well-defined. The main result of \cite{RWC} states (in the notation of the previous section):

\begin{theorem}\label{mainth} In the group ${\rm Aut}(B)$, we have a factorization $$T_{i_1}\circ\ldots\circ T_{i_r}=\prod^\leftarrow_{{\mu\in{\bf Q}}}T_\mu,$$ 
where $$T_\mu(x^d)=x^d\cdot Q_\mu^{\{\_,d\}}(x).$$
Here $Q_\mu^\eta(x)$ denotes the specialization of the series $Q_\mu^\eta(t)$ of Section \ref{section4} from the variables $t_i$ to the variables $x_i$.
\end{theorem}

Let $\Phi\in{\rm Aut}(\Lambda)$ be the map induced on dimension vectors by the inverse Auslander-Reiten translation; $\Phi$ is a Coxeter element of the corresponding Weyl group determined by the property
$$\langle \Phi(d),e\rangle=-\langle e,d\rangle.$$
Then we have $$\{\_,d\}=\langle-({\rm id}+\Phi)d,\_\rangle$$ and thus using Corollary \ref{corsd}:

\begin{corollary} The automorphisms $T_\mu$ of Theorem \ref{mainth} can be written as
$$T_\mu(x_d)=x^d\cdot S_\mu^{-({\rm id}+\Phi)d}(x).$$
\end{corollary}

We specialize to the generalized Kronecker quiver $K_m$ with set of vertices $I=\{i,j\}$ and $m$ arrows from $j$ to $i$. Choose the generators $x=-x_i$ and $y=-x_j$ of $B$; then $B_m={\bf Q}[[x,y]]$ with Poisson bracket $\{x,y\}=mxy$. For $a,b\in{\bf Z}$ with $a,b\geq 0$ and $a+b\geq 1$, we define a Poisson automorphism $T_{a,b}^{(m)}$ of $B$ by
$$T_{a,b}^{(m)}:\left\{\begin{array}{lll}x&\mapsto&x(1-(-1)^{mab}x^ay^b)^{-mb},\\ y&\mapsto&y(1-(-1)^{mab}x^ay^b)^{ma}\end{array}\right.$$
as in \cite[1.4]{KS}. More generally, for an arbitrary series $F(t)\in{\bf Z}[[t]]$ with $F(0)=1$, we define as in \cite[0.1]{GPS}:
$$T_{a,b,F(t)}^{(m)}:\left\{\begin{array}{lll}x&\mapsto&xF(x^ay^b)^{-mb},\\
y&\mapsto&yF(x^ay^b)^{ma}.\end{array}\right.$$
Note that the automorphisms $T_{a,b}^{(m)}$ for fixed slope $a/b$ commute, thus
\begin{equation}\label{commute} \prod_{i\geq 1}(T_{ia,ib}^{(m)})^{d_i}=T_{a,b,F(t)}^{(m)}\end{equation}
for
$$F(t)=\prod_{i\geq 1}(1-((-1)^{mab}t)^i)^{id_i}.$$
We can now use our main results Theorem \ref{genduality}, Theorem \ref{mainth} to confirm \cite[Conjecture 1]{KS}:

\begin{theorem}\label{conj1} Writing
$$T_{1,0}^{(m)}T_{0,1}^{(m)}=\prod^{\leftarrow}_{b/a\mbox{ decreasing}}(T_{a,b}^{(m)})^{d(a,b,m)},$$
we have $d(a,b,m)\in{\bf Z}$ for all $a,b,m$.
\end{theorem}

\proof We choose the stability $\Theta=j^*$ (in fact, the only non-trivial stability, see \cite[5.1]{Rmoduli}).
 By Theorem \ref{mainth}, we have a factorization
\begin{equation}\label{factkron}T_{1,0}^{(m)}T_{0,1}^{(m)}=T_iT_j=\prod^\leftarrow_{{\mu\in{\bf Q}}}T_\mu,\end{equation}
where $$T_\mu(x^d)=x^d\cdot Q_\mu^{\{\_,d\}}(x).$$
Given $\mu\in{\bf Q}$, we write $\mu=b/(a+b)$ for coprime nonnegative $a,b\in{\bf Z}$ and choose integers $c$ and $d$ such that $ac+bd=1$. We have $\Lambda_\mu^+={\bf N}t^{(a,b)}$. Defining
$$G_\mu(t)=Q_\mu^i(t)^cQ_\mu^j(t)^d\in{\bf Z}[[\Lambda_\mu^+]],$$
the proof of \cite[Theorem 6.1]{RWC} shows that
$$G_\mu(t)^a=Q_\mu^i(t)\mbox{ and }G_\mu(t)^b=Q_\mu^j(t).$$
Similarly to Corollary \ref{corsd}, we can find a functional equation for $G_\mu(t)$. We denote $\chi_\mu(k)=\chi(M_{(ka,kb)}^{st}(K_m))$ for $k\geq 1$ and $N=-\langle (a,b),(a,b)\rangle=mab-a^2-b^2$ and apply the first formula of Theorem \ref{t42}:
\begin{eqnarray*}
G_\mu(t)&=&Q_\mu^i(t)^cQ_\mu^j(t)^d=Q_\mu^{(c,d)}(t)\\
&=&\prod_{k\geq 1}R^{(ka,kb)}(t)^{\chi_\mu(k)\cdot k\cdot(ac+bd)}\\
&=&\prod_{k\geq 1}R^{(ka,kb)}(t)^{k\chi_\mu(k)}.\end{eqnarray*}
Applying the second formula of Theorem \ref{t42}, this yields
\begin{eqnarray*}
G_\mu(t)&=&\prod_{k\geq 1}(1-t^{(ka,kb)}\prod_{l\geq 1}R^{(la,lb)}(t)^{klN\chi\mu(l)})^{-l\chi_\mu(l)}\\
&=&\prod_{l\geq 1}(1-t^{(ka,kb)}G_\mu(t)^{kN})^{-k\chi_\mu(k)}.\end{eqnarray*}
Thus, the series $G_\mu(t)$ fulfills the functional equation
\begin{equation}\label{777}G_\mu(t)=\prod_{k\geq 1}(1-(t^{(a,b)}G_\mu(t)^N)^k)^{-k\chi_\mu(k)}.\end{equation}
By Theorem \ref{genduality}, $G_\mu(t)$ admits a factorization
\begin{equation}\label{888}G_\mu(t)=\prod_{k\geq 1}(1-((-1)^Nt^{(a,b)})^k)^{kd_\mu(k)}\end{equation}
for $d_\mu(k)\in{\bf Z}$ for all $k\geq 1$.\\[1ex]
Defining $F_\mu(t)\in{\bf Z}[[t]]$ by $F_\mu((-1)^{a+b}t^{(a,b)})=G_\mu(t)$, we have
\begin{equation}\label{fact16}T_\mu=T_{a,b,F_\mu(t)}^{(m)}\end{equation}
(the sign appearing due to the convention $x=-x_i$, $y=-x_j$) and
\begin{eqnarray*}
F_\mu(t)&=&\prod_{k\geq 1}(1-((-1)^{N+a+b}t)^k)^{kd_\mu(k)}\\
&=&\prod_{k\geq 1}(1-((-1)^{mab}t)^k)^{kd_\mu(k)}.\end{eqnarray*}
By (\ref{commute}) and (\ref{fact16}), this yields
$$T_\mu=\prod_{k\geq 1}(T_{ka,kb}^{(m)})^{d_\mu(k)}.$$
Together with the factorization (\ref{factkron}), this yields the factorization claimed in the theorem, with $d(ka,kb,m)=d_\mu(k)$.\hb

Using a result result of T. Weist, we can also confirm a conjecture in \cite[1.4]{KS} concerning the diagonal term of the factorization in Theorem \ref{conj1}:

\begin{theorem} For all $k\geq 1$, we have
$$d(k,k,m)=\frac{1}{(m-2)k^2}\sum_{i|k}\mu(k/i)(-1)^{mi+1}\binom{(m-1)^2i-1}{i}.$$
\end{theorem}

\proof By \cite[6.2]{Weist}, we have $\chi(M_{d,d}^{st}(K_m))=0$ for $d\geq 2$, whereas $M_{1,1}^{st}(K_m)\simeq{\bf P}^{m-1}$. In the notation of (\ref{777}), (\ref{888}) above, we can apply the example at the end of the previous section with $b=-m$ and $N=m-2$ and arrive at the claimed formula.\hb

\end{document}